\documentclass[10pt]{article}

\usepackage{t1enc}
\usepackage{amsmath,amssymb,amsopn}
\usepackage{theorem}

\theorembodyfont{\sl}
\newtheorem{theorem}{Theorem}[section]
\newtheorem{definition}[theorem]{Definition}
\newtheorem{lemma}[theorem]{Lemma}

\theorembodyfont{\rm}\newtheorem{remark}[theorem]{Remark}
\theorembodyfont{\rm}

\begin{document}

\newcommand{\rot}{\mathop{\textrm{rot}}\nolimits}
\renewcommand{\div}{\mathop{\textrm{div}}\nolimits}
\newcommand{\grad}{\mathop{\textrm{grad}}\nolimits}
\renewcommand{\Re}{\mathop{\textrm{Re}}\nolimits}
\renewcommand{\Im}{\mathop{\textrm{Im}}\nolimits}
\newcommand{\rank}{\mathop{\textrm{rank}}\nolimits}
\renewcommand{\d}{\mathsf{d}}
\newcommand{\hodge}{\mathop{\star}\nolimits}

\title{Connections between optimal constants in some norm inequalities for differential forms}
\author{S{\'a}ndor Zsupp{\'a}n\\ \small{Berzsenyi D{\'a}niel Evang{\'e}likus Gimn{\'a}zium}\\
\small{H-9400, Sopron, Sz{\'e}chenyi t{\'e}r 11., Hungary}\\
\small{zsuppans@gmail.com}}
\maketitle

\begin{abstract}
We derive an improved Poincar{\'e} inequality in connection with the Babu\v{s}ka-Aziz and Friedrichs-Velte inequalities for differential forms by estimating the domain specific optimal constants figuring in the respective inequalities with each other provided the domain supports the Hardy inequality.
We also derive upper estimates for the constants of a star-shaped domain by generalizing the known Horgan-Payne type estimates for planar and spatial domains to higher dimensional ones.
\end{abstract}



\section{Introduction}

The improved Poincar{\'e} inequality for the gradient developed in \cite{hurri1994} has been recently proved in \cite{jiang2014} to be equivalent to the solvability of the divergence equation, that is, to the validity of the Babu\v{s}ka-Aziz inequality for the divergence on a general class of domains.
If we consider the $L_2$-norm, then there is also a simple connection between the domain specific optimal constants in the respective inequalities provided the problem domain supports the Hardy inequality, see \cite{duran2012}.
In \cite{zsuppan2017} with the same method as in the divergence-gradient case an analogous improved Poincar{\'e} inequality for the rotation was derived in connection with the Babu\v{s}ka-Aziz inequality for the rotation and with the corresponding Friedrichs-Velte inequality on spatial domains.
On the other hand equivalence between the Babu\v{s}ka-Aziz and Friedrichs-Velte inequalities has been proved in \cite{costabeldauge2015} for planar domains in the divergence-gradient case and it has been generalized later in \cite{costabel2015} for differential forms on arbitrary dimensional domains.

In this paper we derive relations between the optimal domain specific constants between the Friedrichs-Velte inequality and an improved Poincar{\'e} inequality for differential forms using the framework of \cite{costabel2015}.
Herewith we obtain a generalization of the improved Poincar{\'e} inequality for the case of differential forms, that is, using the exterior derivative instead of the gradient.
We also develop a generalization of the Horgan-Payne type estimates of the Friedrichs-Velte constant of a star-shaped planar \cite{costabeldauge2015,horganpayne1983} or spatial \cite{payne2007} domain for star-shaped domains of arbitrary dimension.

\section{Notation and preliminaries}\label{sec:prelim}

In this paper we use the notations and results of \cite{costabel2015}, which we describe in the following.
Let $\Omega$ be a bounded domain in $\mathbb{R}^n$ ($n\ge 2$) and let $\rho_{\Omega}$ be the distance to the boundary function on $\Omega$, i.e. $\rho_{\Omega}(x)=\text{dist}(x,\partial\Omega)$ for $x\in\Omega$.

For $0\le\ell\le n$ we denote by $\Lambda^{\ell}$ the exterior algebra of $\mathbb{R}^n$ and by $|\cdot|$ the Euclidean norm on $\Lambda^{\ell}$.
$L_2(\Omega,\Lambda^{\ell})$ is the space of differential forms of order $\ell$ with square integrable coefficients.
The norm of $u\in L_2(\Omega,\Lambda^{\ell})$ is defined by
\begin{equation}\label{eq:L2norm}
\|u\|^2=\int_{\Omega}|u|^2
\end{equation}
with the corresponding inner product $\left<\cdot,\cdot\right>$.
In the space $C_0^{\infty}(\Omega,\Lambda^{\ell})$ of smooth differential forms with compact support in $\Omega$ we denote by $\d$ the exterior derivative and by $\d^{\ast}$ the coderivative which is the formal adjoint of $\d$ with respect to the $L_2$ scalar product, i.e.
\begin{equation}\label{eq:exteriordiffs}
\d: C_0^{\infty}(\Omega,\Lambda^{\ell})\to C_0^{\infty}(\Omega,\Lambda^{\ell+1})
\text{ and }
\d^{\ast}: C_0^{\infty}(\Omega,\Lambda^{\ell})\to C_0^{\infty}(\Omega,\Lambda^{\ell-1}),
\end{equation}
where we also have $\d=0$ for $\ell=n$ and $\d^{\ast}=0$ for $\ell=0$.
Let $H_0^1(\Omega,\Lambda^{\ell})$ be the completion of the space $C_0^{\infty}(\Omega,\Lambda^{\ell})$ under the semi-norm
\begin{equation}\label{eq:H1seminorm}
|u|_1=\left<\Delta u,u\right>^{\frac{1}{2}}
\end{equation}
defined by the Hodge-Laplacian
\begin{equation}\label{eq:hodgelaplacian}
\Delta=\d^{\ast}\d+\d\d^{\ast}.
\end{equation}
Denoting by $H^{-1}(\Omega,\Lambda^{\ell})$ the dual space of $H_0^1(\Omega,\Lambda^{\ell})$ with the dual norm
\begin{equation}\label{eqn:dualnorm}
\left|u\right|_{-1}=\sup_{w\in H_0^1(\Omega,\Lambda^{\ell})}\frac{\left<u,w\right>}{\left|w\right|_1},
\end{equation}
we have the following extensions
\begin{eqnarray}
\label{eq:downdiffs}\underline{\d}:H_0^1(\Omega,\Lambda^{\ell})\to L_2(\Omega,\Lambda^{\ell+1})& \text{ and } &\underline{\d}^{\ast}:L_2(\Omega,\Lambda^{\ell+1})\to H^{-1}(\Omega,\Lambda^{\ell}),\\
\label{eq:updiffs}\overline{\d}^{\ast}:H_0^1(\Omega,\Lambda^{\ell})\to L_2(\Omega,\Lambda^{\ell-1}) & \text{ and } & \overline{\d}:L_2(\Omega,\Lambda^{\ell-1})\to H^{-1}(\Omega,\Lambda^{\ell}),
\end{eqnarray}
which we denote simply again by $\d$ and $\d^{\ast}$ below.

We also need the orthogonal complements of the kernels of $\d$ and $\d^{\ast}$:
\begin{eqnarray}
\label{eq:kerupdiffoc}M&=&\left\{u\in L_2(\Omega,\Lambda^{\ell-1})\mid \forall v\in L_2(\Omega,\Lambda^{\ell-1}): \overline{\d}v=0\Rightarrow\left<u,v\right>=0\right\}\\
\label{eq:kerdowncodiffoc}M^{\ast}&=&\left\{u\in L_2(\Omega,\Lambda^{\ell+1})\mid \forall v\in L_2(\Omega,\Lambda^{\ell+1}): \underline{\d}^{\ast}v=0\Rightarrow\left<u,v\right>=0\right\}
\end{eqnarray}
that is $M=\left(\ker\overline{\d}\right)^{\bot}$ and $M^{\ast}=\left(\ker\underline{\d}^{\ast}\right)^{\bot}$.
By Lemma 1.1 in \cite{costabel2015} the operator $\d^{\ast}$ in \eqref{eq:updiffs} maps $H_0^1(\Omega,\Lambda^{\ell})$ to a dense subspace of $M$ and
\begin{equation}\label{eq:codiffM}
u\in M\Rightarrow \d^{\ast}u=0.
\end{equation}
The image of $\d^{\ast}$ coincides with $M$ iff the Babu\v{s}ka-Aziz inequality is satisfied.
\begin{definition}\label{def:BAineq}
The domain $\Omega$ supports the Babu\v{s}ka-Aziz inequality of order $\ell$ if there is a finite positive constant $C_{\ell}$ depending only on $\Omega$ and $\ell$ such that for every $u\in M$ there exists a $w\in H_0^1(\Omega,\Lambda^{\ell})$ such that 
\begin{equation}\label{ineq:BA}
\d^{\ast}w=u\text{ and }|w|_1^2\le C_{\ell}\|u\|^2.
\end{equation}
The least possible constant in \eqref{ineq:BA} is called the Babu\v{s}ka-Aziz constant of order $\ell$ of $\Omega$ and it is denoted by $C_{\Omega,\ell}$.
\end{definition}
The Babu\v{s}ka-Aziz inequality \eqref{ineq:BA} can be equivalently formulated as
\begin{equation}\label{ineq:Lions}
\forall u\in M: \|u\|^2\le C_{\ell,\Omega}|\d u|_{-1}^2,
\end{equation}
or as the inf-sup condition
\begin{equation}\label{eq:infsup}
\inf_{u\in M}\sup_{w\in H_0^1(\Omega,\Lambda^{\ell})}\frac{\left<u,\d^{\ast}w\right>}{\|u\|\cdot |w|_1}=\frac{1}{C_{\ell,\Omega}}.
\end{equation}
This is on the other hand equivalent to
\begin{equation}\label{eq:infsupschur}
\inf_{u\in M}\frac{\left<u,\d^{\ast}\Delta^{-1}\d u\right>}{\left<u,u\right>}=\frac{1}{C_{\ell,\Omega}},
\end{equation}
where the operator 
\begin{equation}\label{def:schurcomplement}
\mathcal{S}=\d^{\ast}\Delta^{-1}\d: L_2(\Omega,\Lambda^{\ell-1})\to L_2(\Omega,\Lambda^{\ell-1})
\end{equation}
is called the Schur complement operator of the Stokes equation, i.e. $\mathcal{S}u=\d^{\ast}w$, $w\in H_0^1(\Omega,\Lambda^{\ell})$ being the solution of $\Delta w=\d u$.
Equation \eqref{eq:infsupschur} means that the least eigenvalue of the restriction of $\mathcal{S}$ onto the subspace $M$ is positive, that is, $\mathcal{S}:M\to M$ is invertible.

According to \cite{costabel2015} the Babu\v{s}ka-Aziz inequality is related to a generalization of the Friedrichs-Velte inequality for conjugate harmonic functions \cite{velte1998}.
The differential forms $u\in L_2(\Omega,\Lambda^{\ell-1})$ and $v\in L_2(\Omega,\Lambda^{\ell+1})$ are called conjugate if they satisfy
\begin{equation}\label{eq:conjugatediffforms}
\d u=\d^{\ast}v.
\end{equation}

\begin{definition}\label{def:FVineq}
The domain $\Omega$ supports the Friedrichs-Velte inequality of order $\ell$ if there is a finite positive constant $\Gamma_{\ell}$ depending only on $\Omega$ and $\ell$ such that for every $u\in M$ and $v\in L_2(\Omega,\Lambda^{\ell+1})$ conjugate in the sense of \eqref{eq:conjugatediffforms} there follows
\begin{equation}\label{ineq:FV}
\|u\|^2\le\Gamma_{\ell}\|v\|^2.
\end{equation}
The least possible constant in \eqref{ineq:FV} is called the Friedrichs-Velte constant of order $\ell$ of $\Omega$ and it is denoted by $\Gamma_{\Omega,\ell}$.
\end{definition}
In inequality \eqref{ineq:FV} we can assume $v\in M^{\ast}$ because for a given $u\in M$ the element $v\in L_2(\Omega,\Lambda^{\ell+1})$ with minimal $L_2$-norm that satisfies \eqref{eq:conjugatediffforms} also satisfies $v\in M^{\ast}$. Hence similar to \eqref{eq:codiffM} we have for such a differential form
\begin{equation}\label{eq:diffcoM}
v\in M^{\ast}\Rightarrow \d v=0.
\end{equation}
Thus the conjugate differential forms \eqref{eq:conjugatediffforms} with $u\in M$ and $v\in M^{\ast}$ satisfy the equations
\begin{equation}\label{eq:conjdiffeqs}
\d^{\ast}u=0,\quad \d u=\d^{\ast}v,\quad \d v=0,
\end{equation}
which mean that both $u$ and $v$ belong to the kernel of the Laplacian \eqref{eq:hodgelaplacian}, that is, both are harmonic.

The following relation between the Babu\v{s}ka-Aziz and Friedrichs-Velte constants has been proved in \cite{costabel2015}.
\begin{theorem}[\cite{costabel2015}, Theorem 2.1]\label{thm:CoBAFV}
For any bounded open set $\Omega\subset\mathbb{R}^n$ and any $1\le\ell\le n-1$, the Babu\v{s}ka-Aziz constant $C_{\Omega,\ell}$ is finite if and only if the Friedrichs-Velte constant $\Gamma_{\Omega,\ell}$ is finite, and there holds
\begin{equation}\label{eq:CoBAFV}
C_{\Omega,\ell}=\Gamma_{\Omega,\ell}+1.
\end{equation}
\end{theorem}

\section{Main results}\label{sec:main}

\subsection{Improved Poincar{\'e} inequality for differential forms}\label{ssec:diffforms}

In order to formulate the main results let us first define the improved Poincar{\'e} inequality for differential forms.
   
\begin{definition}\label{def:iP}
The domain $\Omega\subset\mathbb{R}^n$ supports the improved Poincar{\'e} inequality of order $\ell$ if there is a finite constant $P_{\ell}$ depending only on $\Omega$ and $\ell$ such that for every $u\in M$ there holds
\begin{equation}\label{ineq:iP}
\|u\|^2\le P_{\ell}\|\rho_{\Omega}\d u\|^2.
\end{equation}
The least possible constant $P_{\Omega,\ell}$ in \eqref{ineq:iP} is called the improved Poincar{\'e} constant of order $\ell$ of $\Omega$.
\end{definition}
This definition is a generalization of the improved Poincar{\'e} inequalities introduced in \cite{hurri1994} because we use the exterior derivative instead of the gradient, however, it is only special case because we only consider the $L_2$-space.

Rephrasing \eqref{ineq:iP} yields
\begin{equation*}
\left<u,u\right>\le P_{\ell}\left<\rho_{\Omega}^2\d u,\d u\right>,
\end{equation*}
where $u\in M$, $\rho_{\Omega}^2\d u\in H_0^1(\Omega,\Lambda^{\ell})$ and $\d u\in H^{-1}(\Omega,\Lambda^{\ell})$.
From the previous inequality there follows
\begin{equation*}
\frac{1}{P_{\ell}}\le\frac{\left<\d^{\ast}\rho_{\Omega}^2\d u, u\right>}{\left<u,u\right>}
\end{equation*}
for each $u\in M$, i.e.
\begin{equation}\label{eq:infsupmodifiedschur}
\frac{1}{P_{\Omega,\ell}}=\inf_{u\in M}\frac{\left<\d^{\ast}\rho_{\Omega}^2\d u, u\right>}{\left<u,u\right>}.
\end{equation}
This equation is similar to \eqref{eq:infsupschur} but now we have the diffusion type operator
\begin{equation}\label{eq:modifiedschurcomplement}
\d^{\ast}\rho_{\Omega}^2\d: M\to M
\end{equation}
instead of the Schur complement \eqref{def:schurcomplement}.
Hence, the improved Poincar{\'e} constant in \eqref{eq:infsupmodifiedschur} is connected to the operator \eqref{eq:modifiedschurcomplement} in the same way as the Babu\v{s}ka-Aziz constant in \eqref{eq:infsupschur} to the Schur complement \eqref{def:schurcomplement}. 
Moreover, the inverse of the multiplication operator with factor $\rho_{\Omega}^{-2}$ figuring in \eqref{eq:modifiedschurcomplement} is connected to the Laplacian via the Hardy inequality
\begin{equation}\label{ineq:Hardy_original}
\int_{\Omega}\frac{1}{\rho_{\Omega}^2}|v|^2\le H_{\Omega}\int_{\Omega}\left|\nabla v\right|^2
\text{ for every }
v\in H_0^1(\Omega),
\end{equation}
which was utilized in \cite{duran2012} for the estimation of the improved Poincar{\'e} constant for the gradient from above by the Babu\v{s}ka-Aziz constant for the divergence. 
Regarding inequality \eqref{ineq:Hardy_original} and various estimations for the Hardy constant we refer to \cite{BEL2015} and the references therein.
In this paper we use the following straightforward generalization of \eqref{ineq:Hardy_original}.

\begin{lemma}\label{lem:hardy_for_diff_forms}
Let the domain $\Omega\subset\mathbb{R}^n$ satisfy the Hardy inequality \eqref{ineq:Hardy_original}.
If $v\in H_0^1(\Omega,\Lambda^{\ell})$, then
\begin{equation}\label{ineq:Hardy_norm}
\left\|\frac{1}{\rho_{\Omega}}v\right\|^2\le H_{\Omega}\left|v\right|_1^2,
\end{equation}
where $H_{\Omega}$ is the Hardy constant of $\Omega$.
\end{lemma}
\textbf{Proof.}
Let be $v(x)=\sum v_{i_1,\dots,i_{\ell}}(x)\d x_{i_1}\wedge\dots\wedge\d x_{i_{\ell}}$, where we have for the components $v_{i_1,\dots,i_{\ell}}\in H_0^1(\Omega)$ and the summation extends over the set of all $\ell$-tuples $(i_1,\dots,i_{\ell})$ with $1\le i_1<\dots<i_{\ell}\le n$.
Use the Hardy inequality \eqref{ineq:Hardy_original} for each component:
\begin{equation}
\int_{\Omega}\frac{1}{\rho_{\Omega}^2}|v_{i_1,\dots,i_{\ell}}|^2\le H_{\Omega}\int_{\Omega}|\nabla v_{i_1,\dots,i_{\ell}}|^2.
\end{equation}
In view of
\begin{equation}
\left\|\frac{1}{\rho_{\Omega}}v\right\|^2=\sum\int_{\Omega}\frac{1}{\rho_{\Omega}^2}|v_{i_1,\dots,i_{\ell}}|^2
\text{ and }
|v|_1^2=\sum\int_{\Omega}|\nabla v_{i_1,\dots,i_{\ell}}|^2,
\end{equation}
the assertion of Lemma \ref{lem:hardy_for_diff_forms} follows by summation for all indices.
\hfill$\Box$

\begin{lemma}\label{lem:BAFV->iP}
If the domain $\Omega\subset\mathbb{R}^n$ supports the Babu\v{s}ka-Aziz and Hardy inequalities \eqref{ineq:BA} and \eqref{ineq:Hardy_norm}, respectively, then $\Omega$ also supports the improved Poincar{\'e} inequality \eqref{ineq:iP}.
Moreover, we have
\begin{equation}\label{ineq:P<HC}
P_{\Omega,\ell}\le H_{\Omega}C_{\Omega,\ell}
\end{equation}
between the optimal constants in the respective inequalities.
\end{lemma}
\textbf{Proof.}
The proof is essentially the same as that of Theorem 5.3 in \cite{duran2012} but we repeat it here for the convenience of the reader.
Let be $u\in M$.
According to \cite{costabel2015} by the Babu\v{s}ka-Aziz inequality there exists a differential form $w\in H_0^1(\Omega,\Lambda^{\ell})$ such that $\d^{\ast}w=u$ and $|w|_1^2\le C_{\Omega,\ell}\|u\|^2$.
\begin{eqnarray*}
\|u\|^2&=&\left<u,\d^{\ast}w\right>=\left<\d u,w\right>\le
\left\|\frac{1}{\rho_{\Omega}}w\right\|\cdot\|\rho_{\Omega}\d u\|\le H_{\Omega}^{\frac{1}{2}}|w|_1\|\rho_{\Omega}\d u\|\\
&\le &H_{\Omega}^{\frac{1}{2}}C_{\Omega,\ell}^{\frac{1}{2}}\|u\|\cdot\|\rho_{\Omega}\d u\|
\end{eqnarray*}
There follows 
\begin{equation}
\|u\|^2\le H_{\Omega}C_{\Omega,\ell}\|\rho_{\Omega}\d u\|^2
\end{equation}
for each $u\in M$, which implies that $\Omega$ supports the improved Poincar{\'e} inequality \eqref{ineq:iP}.
In addition we obtain \eqref{ineq:P<HC} and using \eqref{eq:CoBAFV} there also follows
$P_{\Omega,\ell}\le H_{\Omega}\left(1+\Gamma_{\Omega,\ell}\right)$.
\hfill$\Box$

For the opposite direction we adapt the proof of Lemma 3.1 from the paper \cite{zsuppan2017} of the present author.

\begin{lemma}\label{lem:iP->BAFV}
If the domain $\Omega\subset\mathbb{R}^n$ supports the improved Poincar{\'e} inequality \eqref{ineq:iP}, then $\Omega$ also supports the Friedrichs-Velte inequality \eqref{ineq:FV}.
Moreover, we have
\begin{equation}\label{ineq:Gamma<4P}
\Gamma_{\Omega,\ell}\le 4P_{\Omega,\ell}
\end{equation}
between the optimal constants in the respective inequalities.
\end{lemma}
\textbf{Proof.}
Let $u\in M$ and $v\in M^{\ast}$ be conjugate in the sense of \eqref{eq:conjugatediffforms}.
\begin{eqnarray*}
\|\rho_{\Omega}\d u\|^2&=&\left<\rho_{\Omega}^2\d u,\d u\right>=\left<\rho_{\Omega}^2\d u,\d^{\ast} v\right>=\left<\d\left(\rho_{\Omega}^2\d u\right), v\right>\\
&=&\left<\d\left(\rho_{\Omega}^2\right)\wedge\d u+\rho_{\Omega}^2\d\d u,v\right>=2\left<\rho_{\Omega}\d\rho_{\Omega}\wedge\d u,v\right>\\
&=&2\left<\d\rho_{\Omega}\lrcorner v,\rho_{\Omega}\d u\right>\le 2\|\d\rho_{\Omega}\lrcorner v\|\cdot\|\rho_{\Omega}\d u\|,
\end{eqnarray*}
where $\d\rho_{\Omega}\lrcorner v$ is the interior product (contraction) of the differential $\ell+1$ form $v$ with the vector identified with the 1-form $\d\rho_{\Omega}$, i.e.
\begin{equation*}
\d\rho_{\Omega}\lrcorner v=\sum v_{i_1,\dots,i_{\ell+1}}\sum_{k=1}^{\ell+1}(-1)^{k-1}\left(\partial_{i_k}\rho_{\Omega}\right)\d x_{i_1}\wedge\dots\wedge\widehat{\d x_{i_k}}\wedge\dots\wedge\d x_{i_{\ell+1}},
\end{equation*}
for the $\ell+1$ form $v=\sum v_{i_1,\dots,i_{\ell+1}}\d x_{i_1}\wedge\dots\wedge\d x_{i_{\ell+1}}$, where $\widehat{\d x_{i_k}}$ denotes the omission of the corresponding factor and the first summation ranges over all possible $\ell+1$-tuples of indices $1\le i_1<\dots<i_{\ell+1}\le n$.
We obtain
\begin{equation*}
\|\rho_{\Omega}\d u\|\le 2\|\d\rho_{\Omega}\lrcorner v\|,
\end{equation*}
the right-hand side of which we estimate by the norm of $v$.
The Leibniz rule for the interior product gives
\begin{equation*}
\d\rho_{\Omega}\wedge\left(\d\rho_{\Omega}\lrcorner v\right)=
|\d\rho_{\Omega}|^2v-\d\rho_{\Omega}\lrcorner\left(\d\rho_{\Omega}\wedge v\right),
\end{equation*} 
which implies by identity (2.8) in \cite{costabel2009} that
\begin{equation*}
\|\d\rho_{\Omega}\lrcorner v\|^2=
\left<v,v-\d\rho_{\Omega}\lrcorner\left(\d\rho_{\Omega}\wedge v\right)\right>=
\|v\|^2-\|\d\rho_{\Omega}\wedge v\|^2,
\end{equation*}
where we have also used that $|\d\rho_{\Omega}|^2=|\nabla\rho_{\Omega}|^2=1$ almost everywhere in $\Omega$ by the eikonal equation for the boundary distance function.
There follows
\begin{equation*}
\|\d\rho_{\Omega}\lrcorner v\|\le\|v\|
\end{equation*}
by omitting the positive term $\|\d\rho_{\Omega}\wedge v\|^2$ on the right-hand side.
Hence, if $\Omega$ supports the improved Poincar{\'e} inequality \eqref{ineq:iP} then
\begin{equation}
\|u\|^2\le 4 P_{\Omega,\ell}\|v\|^2,
\end{equation}
that is, $\Omega$ also supports the Friedrichs-Velte inequality \eqref{ineq:FV}, moreover, we obtain \eqref{ineq:Gamma<4P} for the involved constants.
By \eqref{eq:CoBAFV} we also have $C_{\Omega,\ell}\le 1+4 P_{\Omega,\ell}$.
\hfill$\Box$

Lemma \ref{lem:BAFV->iP} and Lemma \ref{lem:iP->BAFV} imply the following
\begin{theorem}\label{thm:iP<->BAFV}
If the bounded domain $\Omega\in\mathbb{R}^n$ supports the Hardy inequality \eqref{ineq:Hardy_norm} then $\Omega$ supports the Friedrichs-Velte inequality \eqref{ineq:FV} and simultaneously the Babu\v{s}ka-Aziz inequality \eqref{ineq:BA} if and only if $\Omega$ supports the improved Poincar{\'e} inequality \eqref{ineq:iP}.
Moreover, we have
\begin{equation}\label{ineq:iP<->BAFV}
\frac{1}{4}\left(C_{\Omega,\ell}-1\right)=\frac{1}{4}\Gamma_{\Omega,\ell}\le 
P_{\Omega,\ell}\le 
H_{\Omega}C_{\Omega,\ell}=H_{\Omega}\left(1+\Gamma_{\Omega,\ell}\right)
\end{equation}
for the domain specific constants in the corresponding inequalities.
\hfill$\Box$
\end{theorem}

According to Proposition 3.1.~of \cite{costabel2015} the Babu\v{s}ka-Aziz and the Friedrichs-Velte constants are finite for bounded Lipschitz domains, the Hardy constant of which is also finite.
Hence Theorem \ref{thm:iP<->BAFV} holds for Lipschitz domains with finite constants, that is, Lipschitz domains support the improved Poincar{\'e} inequality \eqref{ineq:iP} for differential forms.

\begin{remark}\label{rem:case_l=1}
Setting $n\ge 2$ and $\ell=1$ as in the Section 4.1.~of \cite{costabel2015}, we have $\d=\overline{\d}=\grad$, $\d^{\ast}=\overline{\d}^{\ast}=-\div$ and the subspace $M=\left(\ker\grad\right)^{\bot}$ consists of functions with vanishing integral over each connected component of $\Omega$.
In this case, according to \cite{hurri1994, jiang2014}, the improved Poincar{\'e} inequality \eqref{ineq:iP} holds for a larger class of domains including John domains.
Moreover, as proved in \cite{jiang2014}, for simply connected planar domains being a John domain is equivalent with the simultaneous finiteness of the investigated constants.
Theorem \ref{thm:iP<->BAFV} gives the additional information that for such domains
\begin{equation*}
\frac{1}{4}\Gamma_{\Omega,\ell}\le 
P_{\Omega,\ell}\le 
16C_{\Omega,\ell}
\end{equation*}
because $H_{\Omega}\le 16$, see \cite{BEL2015}.
\hfill$\Box$
\end{remark}

\begin{remark}\label{rem:hardyconstant}
As one sees from \eqref{ineq:iP<->BAFV} finiteness of the Hardy constant is needed only for estimating the improved Poincar{\'e} constant with the Babu\v{s}ka-Aziz constant.
Finding exact values of the Hardy constant of domains is not so easy, however, for example in case of convex domains in any dimensions we have $H_{\Omega}=4$.
More accessible are estimations of the Hardy constant in terms of geometric characteristics of the domain, e.g. for $\Omega$ satisfying a $\theta$-cone condition we have
\begin{equation*}
H_{\Omega}\le\frac{16}{n\cdot\omega\left(\frac{\sin\theta}{2}\right)},
\text{ where }
\omega\left(\alpha\right)=\frac{\int_0^{\arcsin\alpha}\sin^{n-2}t\,\mathrm{d}t}{2\int_0^{\frac{\pi}{2}}\sin^{n-2}t\,\mathrm{d}t}.
\end{equation*}
For this and other estimations c.f. \cite{BEL2015} and the references given there.
\hfill$\Box$
\end{remark}

\subsection{Estimations for star-shaped domains}

In this section we give an upper estimation for the Friedrichs-Velte constant $\Gamma_{\Omega,1}$ of a domain $\Omega\subset\mathbb{R}^n$ ($n\ge 2$) star-shaped with respect to a ball $B\subset\Omega$.

Note, that there already exist such upper estimations for planar and spatial star-shaped domains. 
For planar domains the differential form $v\in L_2(\Omega,\Lambda^{2})$ specializes to a scalar valued function as its conjugate pair $u$, and the investigation of the associated conjugate pair $u^2-v^2$ and $2uv$ leads to the estimation of the Friedrichs-Velte constant, c.f. \cite{costabeldauge2015,horganpayne1983}.
On the other hand, for spatial domains $v\in L_2(\Omega,\Lambda^{2})$ can be represented by a vector function, \eqref{eq:conjdiffeqs} becomes $\grad u=\rot v$ and $\div v=0$, but the associated functions $u^2-|v|^2$ and $2uv$ are no longer conjugate harmonic.
They satisfy instead identity (2.9) of \cite{payne2007}, which reads in the notation of the present paper as
\begin{equation}\label{eq:payneequation}
\nabla\left(u^2-|v|^2\right)=\rot\left(2uv\right)-2\div\left(v\otimes v\right).
\end{equation}
This identity constitutes the basis for the estimation of the Friedrichs-Velte constant of a spatial star-shaped domain, c.f. \cite{payne2007}.

In order to derive the intended estimation for $\Gamma_{\Omega,1}$ we need a generalization of the identity \eqref{eq:payneequation} when $u\in L_2(\Omega,\Lambda^{0})$ is further a harmonic function normalized by $\int_{\Omega}u=0$ but $v\in L_2(\Omega,\Lambda^{2})$ is a skew-symmetric second order tensor which can be represented by an antisymmetric matrix of order $n(n-1)/2$.

Setting
\begin{equation}\label{eq:v}
v=\sum_{i<j}v_{ij}\d x_i\wedge\d x_j
\end{equation}
we calculate
\begin{eqnarray}
\label{eq:dastv}
\d^{\ast}v&=&\sum_{i=1}^n\left(\sum_{j=1}^n\partial_j v_{ij}\right)\d x_i,\\
\label{eq:dv}
\d v&=&\sum_{i<j<k}\left(\partial_i v_{jk}-\partial_j v_{ik}+\partial_k v_{ij}\right)\d x_i\wedge\d x_j\wedge\d x_k,
\end{eqnarray}
and the equaltions \eqref{eq:conjdiffeqs} become
\begin{equation}\label{eq:conjdiffeqs_2}
\partial_i u=\sum_{j=1}^n\partial_j v_{ij}\text{ and }\partial_i v_{jk}-\partial_j v_{ik}+\partial_k v_{ij}=0.
\end{equation}
Using the first equation of \eqref{eq:conjdiffeqs_2} we calculate
\begin{equation*}
\begin{split}
\d^{\ast}\left(uv\right)&=\sum_{i=1}^n\left(\sum_{j=1}^n\partial_j \left(uv_{ij}\right)\right)\d x_i=
\sum_{i=1}^n\left(\sum_{j=1}^n u\partial_j v_{ij}+v_{ij}\partial_j u \right)\d x_i\\
&= u\,\d^{\ast}v+\sum_{i=1}^n\left(\sum_{j=1}^n\sum_{k=1}^n v_{ij}\partial_kv_{jk}\right)\d x_i,
\end{split}
\end{equation*}
which leads to
\begin{equation*}
\begin{split}
\d^{\ast}\left(uv\right) &= \frac{1}{2}\d\left(u^2+|v|^2\right)+
\sum_{i=1}^n\sum_{j=1}^n\left(
\sum_{k=1}^n v_{ij}\partial_k v_{jk}-\sum_{k=1}^{j-1} v_{jk}\partial_i v_{jk}
\right)\d x_i \\
&= \sum_{i=1}^n\left(\sum_{k=1}^n
\partial_k\left(\sum_{j=1}^n v_{ij}v_{jk}\right)-\sum_{j=1}^n v_{jk}\partial_k v_{ij}-\sum_{j=k+1}^nv_{jk}\partial_i v_{jk}\right)\d x_i.
\end{split}
\end{equation*}
The first term is the divergence of the single dot product of the second order tensor $v$ with itself, while the other terms can be further evaluated to zero using the second equation of \eqref{eq:conjdiffeqs_2} and the skew-symmetry of $v$:
\begin{multline*}
\sum_{k=1}^n\left(\sum_{j=1}^n v_{jk}\partial_k v_{ij}+\sum_{j=k+1}^nv_{jk}\partial_i v_{jk}\right)
=\sum_{j=1}^n\left(\sum_{k=1}^{j-1}v_{kj}\partial_j v_{ik}+\sum_{k=j+1}^n v_{kj}\partial_i v_{kj}\right)\hfill\\
=\sum_{j=1}^n\sum_{k=1}^{j-1}v_{kj}\partial_j v_{ik}-\sum_{j=1}^n\sum_{k=j+1}^nv_{jk}\left(\partial_j v_{ik}-\partial_i v_{jk}\right)\hfill\\
=\sum_{j=1}^n\sum_{k=1}^{j-1}v_{kj}\partial_j v_{ik}-\sum_{j=1}^n\sum_{k=j+1}^n v_{jk}\partial_kv_{ij}
=\sum_{k=1}^n\sum_{j=1}^{k-1}v_{jk}\partial_k v_{ij}-\sum_{j=1}^n\sum_{k=j+1}^n v_{jk}\partial_kv_{ij}\hfill\\
=\sum_{j=1}^n\sum_{k=j+1}^n v_{jk}\partial_k v_{ij}-\sum_{j=1}^n\sum_{k=j+1}^n v_{jk}\partial_kv_{ij}=0
\hfill
\end{multline*}
for each $i=1,\dots,n$.
Thus, as a generalization of \eqref{eq:payneequation} we have obtained
\begin{equation}\label{eq:gen_payneequation}
\d^{\ast}\left(2uv\right)-\d\left(u^2\right)-\d\left(|v|^2\right)=2\div\left(v\cdot v\right).
\end{equation}

\begin{remark}\label{rem:CRMTeq}
For planar domains \eqref{eq:gen_payneequation} reduces to the Cauchy-Riemann system $\nabla^{\bot}(2uv)=\nabla(u^2-v^2)$ because \eqref{eq:v} has only one component and the right hand side of \eqref{eq:gen_payneequation} simplifies to $-2\nabla(v^2)$.
For spatial domains \eqref{eq:gen_payneequation} reduces to \eqref{eq:payneequation} because $\div\left(v\cdot v\right)$ becomes $\div\left(v\otimes v\right)-\d\left(|v|^2\right)$ by identifying the 2-form $v_{12}d x_1\wedge\d x_2+v_{13}\d x_1\wedge\d x_3+v_{23}\d x_2\wedge\d x_3$ with the vector $(v_{23},-v_{13},v_{12})^{\top}$.
\hfill$\Box$
\end{remark}

Another ingredient for the estimation of $\Gamma_{\Omega,1}$ is the 1-form
\begin{equation}\label{eq:oneform_psi}
\psi=
\begin{cases}
\sum_{i=1}^n\left(\frac{1}{r^n}-\frac{1}{r_0^n}\right)x_i\d x_i&\text{ for }a< r\le r_0,\\
\sum_{i=1}^n\left(\frac{1}{a^n}-\frac{1}{r_0^n}\right)x_i\d x_i&\text{ for }0\le r\le a,
\end{cases}
\end{equation}
which generalizes the vector field defined by equation (2.6) of \cite{payne2007}.
In \eqref{eq:oneform_psi} the function
\begin{equation}\label{eqn:r_0}
r=r_0(\theta)\text{ for }\theta\in\mathbb{S},
\end{equation}
parametrizes the boundary of $\Omega$ which is star-shaped with respect to a ball $B$ centered in the origin with radius $a$.  
Here $\mathbb{S}$ denotes the domain of the angles of the $n$-dimensional spherical coordinate system, that is $0\le\theta_j\le\pi$ for $j=1,\dots,n-2$ and $0\le\theta_{n-1}\le 2\pi$.
%Note that according to Chapter 1.1.8 in \cite{mazya1985} $r_0$ is a Lipschitz function on $\mathbb{S}$.
%Observe, that $\psi$ is continuous in $\Omega$ and vanishes on $\partial\Omega$.
By partial integration we obtain
\begin{eqnarray*}
\label{eq:u2v2_psi}
\left<\d\left(u^2\right)+\d\left(|v|^2\right),\psi\right>&=&
\int_{\Omega}\frac{n}{r_0^n}\left(u^2+|v|^2\right)-\frac{n}{a^n}\int_{B}\left(u^2+|v|^2\right),\\
\label{eq:2uv_psi}
\left<\d^{\ast}\left(2uv\right),\psi\right>&=&
-\int_{\Omega}2u\frac{n}{r_0^{n+1}}x^{\top}V\nabla r_0,\\
\label{eq:divvv_psi}
\left<\div\left(v\cdot v\right),\psi\right>&=&
-\int_{\Omega}\frac{2}{r_0^n}|v|^2-\int_{\Omega}\frac{n}{r_0^{n+1}}x^{\top}V^2\nabla r_0
+\int_{\Omega\setminus B}\frac{2}{r^n}|v|^2\\
&&+\int_{\Omega\setminus B}\frac{n}{r^{n+2}}x^{\top}V^2x+\frac{2}{a^n}\int_B|v|^2,
\end{eqnarray*}
where $V=(v_{ij})_{n\times n}$ denotes the antisymmetric matrix corresponding to $v$.
Combining these equations with \eqref{eq:gen_payneequation} gives
\begin{equation}\label{eq:payneeq_psi}
\begin{split}
\int_{\Omega}\frac{n}{r_0^n}u^2&=
\frac{n}{a^n}\int_B u^2+\frac{n-4}{a^n}\int_B |v|^2-\int_{\Omega\setminus B}\left(\frac{4}{r^n}|v|^2+\frac{2n}{r^{n+2}}x^{\top}V^2x\right)\\
& -\int_{\Omega}2u\frac{n}{r_0^{n+1}}x^{\top}V\nabla r_0+\int_{\Omega}\frac{2n}{r_0^{n+1}}x^{\top}V^2\nabla r_0+
\int_{\Omega}\frac{4-n}{r_0^n}|v|^2
\end{split}
\end{equation}
From this point on we impose on the harmonic function $u$ the normalization $u(0)=0$ instead of $\int_{\Omega}u=0$.
This can be done because
\begin{equation*}
\int_{\Omega}\left(u-u(0)\right)^2=\int_{\Omega}u^2+|\Omega|u(0)^2\ge\int_{\Omega}u^2\text{ if }\int_{\Omega}u=0.
\end{equation*}
These two normalizations are equivalent by the mean value theorem if $\Omega$ itself is a ball.
Similar to equation (2.14) in \cite{payne2007} we obtain using \eqref{eq:conjdiffeqs_2} and $v_{ii}=0$ the estimation
\begin{equation*}
\begin{split}
\Delta\left(u^2-(n-1)|v|^2\right)&=2\left(|\nabla u|^2-(n-1)\sum_{i<j}|\nabla v_{ij}|^2\right)\\
&\le 2\left(\sum_{i=1}^n\left(\sum_{j=1}^n\partial_j v_{ij}\right)^2-(n-1)\sum_{i<j}|\nabla v_{ij}|^2\right)\\
&= 2(n-1)\sum_{i<j}\left(\left(\partial_i v_{ij}\right)^2+\left(\partial_j v_{ij}\right)^2-|\nabla v_{ij}|^2\right)
\le 0,
\end{split}
\end{equation*}
which implies by the mean value theorem and the normalization $u(0)=0$ that $\Gamma_{B,1}\le n-1$ for an $n$-dimensional ball.
In fact $\Gamma_{B,1}=n-1$, c.f. \cite{costabeldauge2015,horganpayne1983,payne2007}.
Substituting this into \eqref{eq:payneeq_psi} there follows
\begin{equation*}
\begin{split}
\int_{\Omega}\frac{n}{r_0^n}u^2&\le
\frac{n^2-4}{a^n}\int_B |v|^2-\int_{\Omega\setminus B}\left(\frac{4}{r^n}|v|^2+\frac{2n}{r^{n+2}}x^{\top}V^2x\right)\\
&-\int_{\Omega}2u\frac{n}{r_0^{n+1}}x^{\top}V\nabla r_0+\int_{\Omega}\frac{2n}{r_0^{n+1}}x^{\top}V^2\nabla r_0+
\int_{\Omega}\frac{4-n}{r_0^n}|v|^2.
\end{split}
\end{equation*}
Using $x^{\top}V^2x=-|Vx|^2$ and $|Vx|^2=|x\lrcorner v|^2\le r^2|v|^2$ we estimate the integral over $\Omega\setminus B$ as
\begin{equation*}
-\int_{\Omega\setminus B}\left(\frac{4}{r^n}|v|^2+\frac{2n}{r^{n+2}}x^{\top}V^2x\right)
%=\int_{\Omega\setminus B}\left(\frac{2n}{r^{n+2}}|Vx|^2-\frac{4}{r^n}|v|^2\right)
\le\int_{\Omega\setminus B}\frac{2n-4}{r^n}|v|^2
\end{equation*}
and we obtain
\begin{equation}\label{ineq:payneineq_psi}
\begin{split}
\int_{\Omega}\frac{n}{r_0^n}u^2&\le
\frac{n^2-4}{a^n}\int_{\Omega} |v|^2+\int_{\Omega}\frac{4-n}{r_0^n}|v|^2\\
&-\int_{\Omega}2u\frac{n}{r_0^{n+1}}x^{\top}V\nabla r_0+\int_{\Omega}\frac{2n}{r_0^{n+1}}x^{\top}V^2\nabla r_0.
\end{split}
\end{equation}

\begin{remark}\label{rem:2D_CHP_bound}
For planar domains $\Omega$ the first term in \eqref{ineq:payneineq_psi} vanishes evidently. The last term also vanishes because $V^2$ is a multiple of the identity matrix and $\nabla r_0$ is orthogonal to $x$ since $r_0$ depends only on the angular variables.
The estimation \eqref{ineq:payneineq_psi} simplifies to
\begin{equation*}\label{ineq:CHP}
\int_{\Omega}\frac{1}{r_0^2}u^2\le
\int_{\Omega}\frac{1}{r_0^2}v^2-\int_{\Omega}2uv\frac{x_1\partial_2r_0-x_2\partial_1r_0}{r_0^3},
\end{equation*}
wherein the estimation of the last term leads to the Horgan-Payne estimation of the Friedrichs constants of $\Omega$, c.f. (4.10) and subsequent equations in \cite{costabeldauge2015} for example.
\hfill$\Box$
\end{remark}
The latter two integrals in \eqref{ineq:payneineq_psi} are estimated as
\begin{eqnarray*}
\int_{\Omega}2u\frac{n}{r_0^{n+1}}x^{\top}V\nabla r_0&\le &
\alpha\int_{\Omega}\frac{n}{r_0^n}u^2+\frac{1}{\alpha}\int_{\Omega}\frac{n}{r_0^n}|v|^2\frac{|r\nabla r_0|^2}{r_0^2},\\
\int_{\Omega}\frac{2n}{r_0^{n+1}}x^{\top}V^2\nabla r_0&=&
-\int_{\Omega}\frac{2n}{r_0^{n+1}}(Vx)^{\top}(V\nabla r_0)\le 
\int_{\Omega}\frac{n}{r_0^n}|v|^2\frac{2|r\nabla r_0|}{r_0}
\end{eqnarray*}
for some parameter $0<\alpha<1$.
Substituting these estimates into \eqref{ineq:payneineq_psi} gives
\begin{equation}\label{ineq:gen_payne}
(1-\alpha)\int_{\Omega}\frac{n}{r_0^n}u^2\le 
\frac{n^2-4}{a^n}\int_{\Omega} |v|^2+\int_{\Omega}\frac{4-n}{r_0^n}|v|^2+
\left(\frac{1}{\alpha}Q^2+2Q\right)\int_{\Omega}\frac{n}{r_0^n}|v|^2,
\end{equation}
where $\frac{|r\nabla r_0|}{r_0}$ depends only on the angular variables of the spherical coordinate system and hence the quantity $Q=\max_{\partial\Omega}\frac{|r\nabla r_0|}{r_0}\ge 0$ depends only on the shape of $\Omega$.
%For example, if $\Omega$ is a ball then $Q=0$.
\begin{remark}\label{rem:3D_P_bound}
In case of three dimensional domains \eqref{ineq:gen_payne} reduces to
\begin{equation*}
\int_{\Omega}u^2\le\frac{2+2Q+\frac{Q^2}{\alpha}}{1-\alpha}\cdot\max\left(\frac{r_0}{a}\right)^3\cdot\int_{\Omega}|v|^2,
\end{equation*}
which leads to
\begin{equation}
\int_{\Omega}u^2\le\left(Q+\sqrt{Q^2+2Q+2}\right)^2\cdot\max\left(\frac{r_0}{a}\right)^3\cdot\int_{\Omega}|v|^2
\end{equation}
by optimizing with respect to $\alpha$.
This constitutes an upper estimation for the Friedrichs-Velte constant of the spatial domain comparable to that of \cite{payne2007}.
\hfill$\Box$
\end{remark}
For any dimensions $n\ge 2$ from \eqref{ineq:gen_payne} there follows
\begin{equation*}
\frac{(1-\alpha)n}{\max_{\partial\Omega}r_0^n}\int_{\Omega}u^2\le
\left(\frac{n^2-4}{a^n}+\frac{4-n}{\max_{\partial\Omega}r_0^n}+\frac{\frac{nQ^2}{\alpha}+2nQ}{a^n}\right)
\int_{\Omega}|v|^2,
\end{equation*}
which implies
\begin{equation*}\label{ineq:gen_payne_eta}
\int_{\Omega}u^2\le 
\eta^n\frac{\frac{n^2-4}{n}-\frac{n-4}{n}\eta^{-n}+\frac{Q^2}{\alpha}+2Q}{1-\alpha}
\int_{\Omega}|v|^2,
\end{equation*}
where $\eta=\max_{\partial\Omega}\frac{r_0}{a}$ is the eccentricity of $\Omega$ with respect to its center of star-shapedness.
By optimizing with respect to $\alpha$ we obtain the upper bound
\begin{equation}\label{ineq:gen_FVBound}
\Gamma_{\Omega,1}\le 
\eta^n\left(
Q+\sqrt{Q^2+2Q+\frac{n^2-4}{n}-\frac{n-4}{n}\eta^{-n}}
\right)^2
\end{equation}
for the Friedrichs-Velte constant.
\begin{remark}\label{rem:FVnDimBall}
For an $n$-dimensional ball $B$ we have $Q=0$ and $\eta=1$ w.r.t. its center.
Substituting these into \eqref{ineq:gen_FVBound} gives the known bound $\Gamma_{B,1}\le n-1$.
\hfill$\Box$
\end{remark}

\begin{remark}\label{rem:nDimBAiPEstim}
The upper estimation \eqref{ineq:gen_FVBound} carries over by Theorem \ref{thm:iP<->BAFV} to the corresponding Babu\v{s}ka-Aziz constant for the divergence and improved Poincar{\'e} constant for the gradient of an $n$-dimensional star-shaped domain.
It improves on the estimation of \cite{galdi1994_1}.
It also generalizes the estimation in \cite{chuawheeden2010} for $P_{\Omega,1}$ of a convex domain for star-shaped ones, however, only in the $L_2$-case, c.f. \cite{chuawheeden2010}.
\hfill$\Box$
\end{remark}

\begin{remark}\label{rem:estim_general_ell}
The estimation \eqref{ineq:gen_FVBound} is valid only for the case $\ell=1$.
In order to obtain an estimation for every $\ell\ge 1$ we have to generalize \eqref{eq:gen_payneequation}.
Instead of $uv$ on the left-hand side of \eqref{eq:gen_payneequation} we could consider $\hodge\left(u\wedge\hodge v\right)$ for a pair $u\in L_2(\Omega,\Lambda^{\ell-1})$ and $v\in L_2(\Omega,\Lambda^{\ell+1})$ conjugate in the sense \eqref{eq:conjdiffeqs}.
%so that the left hand side of \eqref{eq:gen_payneequation} we would be
%\begin{equation*}
%\d^{\ast}\left(2\hodge\left(u\wedge\hodge v\right)\right)-\d\left(|u|^2\right)-\d\left(|v|^2\right).
%\end{equation*}
However, this is beyond the scope of this paper.
\hfill$\Box$
\end{remark}

\subsection*{Concluding remarks}
In this paper, motivated by the results of \cite{costabel2015}, we have derived relations between the domain specific improved Poincar{\'e} constant and the Friedrichs-Velte constant figuring in the corresponding inequalities for differential forms.
This constitutes a generalization of the improved Poincar{\'e} inequalities formulated in \cite{hurri1994} for the gradient, however, only in case of the $L_2$-space on a Lipschitz domain.
We have also developed a generalization of the Horgan-Payne type upper estimations for the Friedrichs-Velte constant of a planar \cite{costabeldauge2015,horganpayne1983} and of a spatial \cite{payne2007} domain for arbitrary dimensional star-shaped domains, which also constitutes a unification of these estimations.
However, this generalization is only valid for conjugate harmonic differential 0- and 2-forms in arbitrary dimensions.
For more general conjugate harmonic differential $\ell-1$- and $\ell+1$-forms we only indicated a possible way of generalization.

\end{document}